\documentclass[12pt]{amsart}
\usepackage{amsmath, color}
\usepackage[all]{xy}
\newif\ifpdf
	\ifx\pdfoutput\undefined
	\pdffalse 
	\else
	\pdfoutput=1 
	\pdftrue
	\fi

	\ifpdf
	\usepackage[pdftex]{graphicx}
	\else
	\usepackage{graphicx}
	\fi

\theoremstyle{plain}
\newtheorem{thm}{Theorem}
\newtheorem{corol}{Corollary}

\newtheorem{lem}{Lemma}
\newtheorem{prop}{Proposition}

\theoremstyle{definition}
\newtheorem{defin}{Definition}
\theoremstyle{remark}
\newtheorem{rem}{Remark}
\theoremstyle{remark}

\def\qed{\hfill\vrule height 5pt width 5 pt depth 0pt}

\begin{document}

\ifpdf
	\DeclareGraphicsExtensions{.pdf, .jpg, .tif}
	\else
	\DeclareGraphicsExtensions{.eps, .jpg}
	\fi

\title{Geometry of Infinitely Generated Veech Groups}
\author{ Pascal Hubert}
\address{Institut de Math\'ematiques de Luminy,
         163 av de Luminy, case 907,
         13288 Marseille cedex 09,
         France} 
\email {hubert@iml.univ-mrs.fr}
\author{Thomas A. Schmidt}
\address{Oregon State University\\ Corvallis, OR 97331}
\email{toms@math.orst.edu}

\date{27 July 2004}


\begin{abstract}  Veech groups uniformize Teichm\"uller geodesic curves in Riemann moduli space.   Recently, examples of infinitely generated Veech groups have been given.  We show that these can even have infinitely many cusps and infinitely many infinite ends.  We further show that  examples exist for which each direction of an infinite end is the limit of directions of inequivalent infinite ends.   

\end{abstract}

\maketitle

\section{Introduction}

  The study of billiards on Euclidean polygonal tables has benefitted greatly from the use of the theory of Riemann surfaces.   A classical unfolding process associates to each rational-angled polygonal table a Riemann surface with holomorphic 1-form, see \cite{KZ}. Using this, Kerkhoff-Masur-Smillie \cite{KMS}  showed that for  almost every direction on a polygonal table, the billiard flow is ergodic.  Soon thereafter, Veech \cite{Vch1} identified a setting where the only non-ergodic directions are those of saddle connections; these examples are identified by an associated group, the Veech group, being a lattice Fuchsian group.     In \cite{Vch4}, Veech asked if in general  these groups can be non-finitely generated.  Only  recently have examples of infinitely generated Veech groups  been given \cite{HS3}, \cite{Mc2}.     Here we study the geometry associated to some of these groups.    
  
  The Riemann moduli space $\mathcal M_g$  is the space of smooth compact Riemann surfaces of genus $g$,  up to biholomorphic equivalence.      The holomorphic 1-forms on a Riemann surface $X$ of genus $g$ form a $g$-dimensional complex vector space, $\Omega(X)$.  Pairs $(X, \omega)$ with $\omega\in \Omega(X)$ form a bundle
 $\Omega \mathcal M_g \to \mathcal M_g$ (of course, $\mathcal M_g$ is not a manifold, but rather an orbifold, thus a fiber is $\Omega(X)/\text{Aut}(X)\,$).   Each pair $(X, \omega)$ (with $\omega$ non-null)  identifies a translation structure  on the real  surface $X$:   local coordinates are defined by integration of $\omega$.    Invoking the real structure of $\mathbb C$, there is a $\text{SL}(2, \mathbb R)$-action on these translation surfaces:  post-composition with the local coordinate functions.    
 
To each $(X, \omega)$, there is thus a  map from $\text{SL}(2, \mathbb R)$ to $\mathcal M_g$, given by sending group elements to corresponding orbit elements.   But,  rotations act so as to leave the  Riemann surface $X$ fixed.    Thus there is an induced a  map from $\text{SO}_{2}(\mathbb R)\backslash \text{SL}(2, \mathbb R)$  to  $\mathcal M_g$.   With respect to the Teichm\"uller metric on $\mathcal M_g$, this map turns out to be  an isometric immersion of the hyperbolic plane $\mathbb H^2 =\text{SL}(2, \mathbb R)/\text{SO}_{2}(\mathbb R)$.   The resulting image is called the  {\em Teichm\"uller complex geodesic} determined by $(X, \omega)$.   We will denote this by $\mathcal G(X, \omega)$.

The $\text{SL}(2, \mathbb R)$-stabilizer of  $(X, \omega)$  induces a subgroup of the oriented isometry group of the hyperbolic plane.   This induced group is the {\em Veech group}, denoted $\text{PSL}(X, \omega)$.   The Teichm\"uller complex geodesic is thus the image of the hyperbolic surface (or possibly orbifold) uniformized by the Veech group;  the map between them is generically injective .    We use this to show:

\begin{thm}\label{thmGeometry} For any genus $g\ge 4$, there exists a complex Teichm\"uller geodesic in the Riemann moduli space  $\mathcal M_g$ that has infinitely many cusps and infinitely many infinite ends.   \end{thm}  

We remind the reader of the definition and basic classification of ends of a hyperbolic surface below.  Informally, an  {\em end}  is a direction in which geodesics can leave compact subsets.    Finite volume hyperbolic surfaces can only have cusps as ends;   briefly, a {\em cusp}, or puncture, corresponds to a conformal punctured unit disk in the surface (which cannot be completed in the surface).   An end of a 
 hyperbolic surface which corresponds to a conformal non-trivial annulus is called a hole, or a flare.   These two types of ends are called finite, an {\em infinite end} is an end which is not finite.

 We prove our main results using the approach of translation surfaces.   In particular,  in terms that we make precise below, whereas cusps always form a discrete set, 
  we give examples such that any direction giving an infinite end is the limit of directions of inequivalent infinite ends.    A translation surface $(Y, \alpha)$ is called a {\em Veech surface} if  $\text{PSL}(Y, \alpha)$ is a lattice in $\text{PSL}(2, \mathbb R)$; the definition of a non-periodic connection point can be found in Subsection \ref{sadParSpecPts} 

\begin{thm}  \label{infFatBoys}  Let $(X, \omega)$ be a translation surface which is a finite ramified covering of a Veech surface $(Y, \alpha)$, with ramification above some non-periodic connection point and possibly the singularities of $(Y, \alpha)$.    Then  $\text{PSL}(X, \omega)\backslash \mathbb H^2$ has infinitely many cusps and infinitely many infinite ends.      Choose a Dirichlet fundamental domain for $\text{PSL}(X, \omega)$ acting upon the hyperbolic plane $\mathbb H^2$;   this determines an infinite set $\Xi$ of points in $\partial \mathbb H^2$, each  representing a distinct infinite end.   The set  $\Xi$ has no isolated points.    
\end{thm}  

 That  $\text{PSL}(X, \omega)$, for $(X, \omega)$ as in Theorem \ref{infFatBoys}, is infinitely generated is the main result of \cite{HS3}.   For the $(X, \omega)$ of Theorem \ref{infFatBoys}, the  complicated geometry at infinity of $\text{PSL}(X, \omega)\backslash \mathbb H^2$ is projected to a small sublocus of the Riemann moduli space.

\begin{thm}\label{thmBigCusp}   Let  $(X,\omega)$ be as in Theorem \ref{infFatBoys}, and suppose 
that $X$ is of genus $g$.    Then, there is a two-complex dimensional sublocus $E\subset \mathcal M_g$ such that all ends of the complex  Teichm\"uller geodesic generated by $(X,\omega)$ are contained in finitely many  ends of $E$.
 \end{thm}  

\subsection{Outline} 
We prove Theorem \ref{thmGeometry} by restricting to the case where  $(X, \omega)$ satisfies the hypotheses of Theorem \ref{infFatBoys} --- an application of the Riemann-Hurwitz formula shows that we may assume that the corresponding $(Y, \alpha)$ be of genus two, we give explicit examples of such $(Y, \alpha)$ which admit non-periodic connection points in \cite{HS3}, it is now known \cite{C}, \cite{Mc1} that the set of such examples is large --- the existence of infinitely many cusps is proven in \S \ref{manyCusps},  the existence of infinitely many infinite ends is proven in \S \ref{manyEnds}, that every direction of an infinite end is the limit of distinct infinite end directions is also proven in \S \ref{manyEnds}.     We prove Theorem \ref{thmBigCusp} in \S \ref{locEnd}.     In \S 2 we give background material on translation surfaces and certain $\text{SL}(2, \mathbb R)$-actions, in \S 3 we discuss projections of horocycles to the moduli space $\mathcal M_g$.   The main results rest upon the elementary tools presented in \S 4.

\subsection{Comments on a Related Setting} 
That a hyperbolic surface uniformized by an infinitely generated Fuchsian group of the first kind admits infinite ends can be deduced  already  from \cite{F}, Theorem 16.  
Thus, the examples of infinitely generated Veech group given by McMullen \cite{Mc2} must also have these.    Each of McMullen's examples has its associated Teichm\"uller complex geodesics lying  densely on a symmetric Hilbert modular surface, 
all ends must lie in the finitely many ends of the Hilbert modular surface.    Thus, if here also there are infinitely many infinite ends, then they may well show the geometry of the type we discuss here. 

\subsection{Thanks} 
It is a pleasure to thank F. Dal'bo, G. Forni, and J. Smillie for helpful discussions.   In particular, we thank G. Forni for his ideas which are reflected in \S \ref{projHoro}.    We thank A. Eskin and others, see \cite{EMM}, for continued interest in these matters.

\section{Background: $\text{SL}(2, \mathbb R)$-Action on Translation Surfaces }

In order to establish notation, we review basic notions.  See say \cite{MT} for more details;  \cite{HS1}, \cite{HS2} and \cite{HS3} use related notation and notions.    See also   \cite{Mc1}, \cite{Mc2}, \cite{Mc3} for recent related developments.   

\subsection{Translation Surfaces and Affine Diffeomorphisms} 

A {\em translation surface}  is a real surface with an atlas such that all transition functions are translations.   As usual, we consider maximal atlases.  Given a Riemann surface $X$ and a holomorphic 1-form $\omega \in \Omega(X)$, integration of $\omega$ induces local coordinates on $X$ away from the zeros of $\omega$; regarding $\mathbb C$ with its standard structure as $\mathbb R^2$, transition functions are easily seen to be translations.  In particular,  directions of flow on a translation surface 
are well-defined, by pull-back from the plane.   One easily checks that the translation surface associated to $(X, \omega)$ is well-defined up to  isomorphism of translation surface (such isomorphism is realized by local pure translations), see Section 2.1  of \cite{HS3} for a more formal discussion.   

The holomorphic 1-form $\omega$ induces a flat metric on $X$ off of its zeros; when the metric is completed, a zero of multiplicity $k$  gives a conical singular point, of cone angle $2(k+1) \pi$.     We consider only translation surfaces of this type; the  notation $(X, \omega)$  hence denotes both  the ordered pair  of Riemann surface $X$ with holomorphic 1-form $\omega$,  and the translation surface structure,  with which $\omega$ equips $X$.   

Given a translation surface $(X, \omega)$, and an element $A \in \text{SL}(2, \mathbb R)$, we can post-compose the coordinate functions of the charts of the atlas of  $(X, \omega)$ by $A$.   It is easily seen that this gives again a translation surface,  denoted by $A \circ (X, \omega)$.

Given a translation surface  $(X, \omega)$, let $X'$ denote the surface which arises from deleting the singularities of  $(X, \omega)$.   An {\em affine diffeomorphism} $f$ of $(X, \omega)$ is a
 homeomorphism of $X$ that restricts to a diffeomorphism of $X'$ whose differential at all points is a constant  element of $\text{SL}(2, \mathbb R)$.   In other terms, $f$ induces a real-affine map in the 
charts of the $\omega$-atlas on $X$,  of constant (Jacobian) derivative.   The orientation-preserving affine diffeomorphisms form a group, denoted $\text{Aff}^{+}(X, \omega)$.     The subgroup of $\text{SL}(2, \mathbb R)$ which arises as the derivatives of affine diffeomorphisms is denoted by
$\text{SL}(X, \omega)$.

  If $\phi, \psi \in \text{Aff}^{+}(X, \omega)$ have the same derivative, then the diffeomorphism $\phi\circ \psi^{-1}$ acts by pure translations.  This diffeomorphism thus preserves the conformal structure on the underlying translation surface; that is, $\phi\circ \psi^{-1}$ is  a conformal automorphism of the Riemann surface $X$ (that preserves $\omega$).   The kernel of the group homomorphism from $\text{Aff}^{+}(X, \omega)$ to $\text{SL}(X, \omega)$ is hence certainly finite whenever $X$ is of genus greater than one.     By pull-back on charts, each $\phi \in  \text{Aff}^{+}(X, \omega)$ also gives  a  translation structure on the underlying surface.  Indeed, if $A$ is the constant derivative of $\phi$ with respect to the $\omega$-charts, then $\phi$ is an isomorphism of the translation surface $A \circ (X, \omega)$ with  $(X,\omega)$.    Thus, the $\text{SL}(2, \mathbb R)$-stabilizer of $(X, \omega)$  is $\text{SL}(X, \omega)$.

\subsection{Teichm\"uller Complex Geodesic, Veech Group} \label{VchGrp} 

Fix some $(X, \omega)$, and let $F: \text{ SL}(2, \mathbb R) \to \Omega \mathcal M_g$ be given by $A \mapsto A \circ (X, \omega)$;   let $\pi:  \Omega \mathcal M_g \to \mathcal M_g$ denote the standard projection.   Note that the underlying Riemann  surface $Y$  is common to the $\text{SO}_{2}(\mathbb R)$-orbit of any $(Y,\alpha)$ --- the image of $(Y, \alpha)$ under a rotation is  $(Y, \zeta \cdot \alpha)$, where $\zeta\in S^1$ is the complex number which induces the rotation.  Thus, $\pi\circ F$  induces a map (from right cosets) $f: \text{SO}_{2}(\mathbb R)\backslash\text{ SL}(2, \mathbb R) \to \mathcal M_g$.     The image of $f$ is $\mathcal G(X, \omega)$, the {\em Teichm\"uller complex geodesic} associated to  $(X, \omega)$.   

Right multiplication gives the natural right group action of  $\text{ SL}(2, \mathbb R)$ on the cosets $\text{SO}_{2}(\mathbb R)\backslash\text{ SL}(2, \mathbb R)$.   The center, 
$\{-I, I\} = \langle-I \rangle$, of 
$\text{ SL}(2, \mathbb R)$ is a subgroup of $\text{SO}_{2}(\mathbb R)$, and one finds that  $\text{PSL}(2, \mathbb R) =  \text{ SL}(2, \mathbb R)/ \langle-I \rangle$ acts faithfully.    The $\text{PSL}(2, \mathbb R)$-stabilizer of the map $f$ is the {\em Veech group} $\text{PSL}(X,\omega)$,  the image in $\text{PSL}(2, \mathbb R)$ of $\text{SL}(X,\omega)$.   This is trivially verified:  $f(\, \text{SO}_{2}(\mathbb R).A\,) = f(\,\text{SO}_{2}(\mathbb R).B\,)$ implies that there is some rotation $R$ such that $B^{-1}RA \in \text{SL}(X,\omega)$.  Thus, the class $\text{SO}_{2}(\mathbb R).A$ is in the $\text{PSL}(X,\omega)$-orbit of $\text{SO}_{2}(\mathbb R).B\,$.   That $f$ is single-valued on each such orbit is obvious.    The map $f$ is generically injective; it can only fail to be injective at some $(Y,\alpha)$ if 
 $\text{Aut}(Y)$ is larger than the automorphism group of Riemann surfaces corresponding to neighboring points of $\mathcal M_g$, see \cite{R} and \cite{P}.   When the Veech group normalizes a true 
 surface, then $f$ gives a  ramified cover normalizing, in the sense of algebraic curves,
$\mathcal G(X, \omega)$. 

    In \cite{Vch1}, Veech proved that $\text{PSL}(X,\omega)$ is a discrete subgroup of $\text{PSL}(2, \mathbb R)$.  That is, it is a Fuchsian group.  The celebrated {\em Veech Dichotomy}, also proven in \cite{Vch1},  states that if $\text{PSL}(X, \omega)$ is a lattice --- that is, if $\text{PSL}(X,\omega)\backslash \mathbb H^2$ has finite hyperbolic area ---, then the 
linear flow in any direction of the surface $(X, \omega)$ is either ergodic or periodic. 
One thus says that a translation surface is a {\em Veech surface} if $\text{PSL}(X, \omega)$ is a lattice.   We also say that  $(X, \omega)$ is an {\em arithmetic surface} if 
$\text{PSL}(X, \omega)$ is arithmetic, in that it is $\text{PSL}(2, \mathbb R)$-commensurable to $\text{PSL}(2, \mathbb Z)$.

\subsection{Teichm\"uller Metric} \label{TeichMet} 
The Teichm\"uller metric on $\mathcal M_g$ is found by giving the Teichm\"uller metric on the (ramified) cover Teichm\"uller space $\mathcal T_g$ of $\mathcal M_g$; the Teichm\"uller modular group (equivalently, the mapping class group) acts isometrically on $\mathcal T_g$, with quotient $\mathcal M_g$.   Thus, the Teichm\"uller metric descends to $\mathcal M_g$. We avoid direct discussion of Teichm\"uller space here, the reader may wish to refer to, say, \cite{G}.

If $r$ is a non-zero real, then the translation surfaces $(X, \omega)$ and $(X, r\omega)$ differ only by scaling.  This scaling and the $\text{SL}(2, \mathbb R)$-action commute; furthermore,  the projection map $\pi: \Omega \mathcal M_g \to \mathcal M_g$ ignores such scaling.   Thus, we can restrict our discussion to $\Omega_1 \mathcal M_g$, the collection of area one $(X, \omega)$.

In informal terms, whereas a conformal map of Riemann surfaces takes infinitesimal circles to circles, a {\em quasi-conformal} map of a Riemann surfaces takes infinitesimal circles to ellipses.    The {\em Beltrami coefficient} of a quasi-conformal homeomorphism $\phi: X\to Y$ is computed by $\bar \partial(f_U \circ \phi)/\partial(f_U \circ \phi)$, where $(U, f_U)$ is a local chart in $Y$, and the symbols $\partial$ and $\bar \partial$ denote $\partial/\partial z$ and $\partial/\partial \bar z$, respectively.

Now, if $(X, \omega) \in \Omega_1 \mathcal M_g$ and 
 $A \in \text{SL}(2, \mathbb R)$,  then the identity map on $X$, from $(X, \omega)$ to $A\circ (X, \omega)$ has Beltrami coefficient 
$\mu_A \,\bar \omega/\omega$, where $\mu_A $ is the complex dilatation of the linear map given by $A$ on $\mathbb C$ with its standard real structure.    But, the map sending elements $A\in \text{SL}(2, \mathbb R)$ to $\mu_A$ is constant on right  $\text{SO}(2, \mathbb R)$-cosets.  This can be shown by direct calculation, or see 
Corollary 1 of \cite{EG}.   Indeed, $A \mapsto \mu_A$ injectively identifies $\text{SO}(2, \mathbb R)\backslash \text{SL}(2, \mathbb R)$ with the complex unit disk, $\Delta$.

    The Teichm\"uller metric on $\text{SO}(2, \mathbb R)\backslash \text{SL}(2, \mathbb R)$ is the pull-back of the hyperbolic metric on $\Delta$.      But, as \cite{EG} show in their Lemma 4.2,  the map sending $\Delta \to \Delta$ by $\mu_A \mapsto \mu_{AB}$ is conformal;  it is thus an oriented isometry.   Hence, the Veech group is the stabilizer in the oriented isometry group of $f: \text{SO}_{2}(\mathbb R)\backslash\text{ SL}(2, \mathbb R) \to \mathcal M_g$.   
      
\bigskip 

It is reasonable to identify $\text{SO}(2, \mathbb R)\backslash \text{SL}(2, \mathbb R)$  with $\mathbb H^2$, as does \cite{Mc1}. In this case the stabilizer of the map
$\mathbb H^2 \to \mathcal M_g$ giving the Teichm\"uller complex geodesic is again the Veech group, but acting on the left by m\"obius transformations conjugated with the action of $z\mapsto -\bar z$.  This conjugating action corresponds to $\tau = \begin{pmatrix} -1&0\\ 0 & 1\end{pmatrix}$, the constant derivative of a coset representative of the index two subgroup $\text{Aff}^{+}(X,\omega)$  in the full group of affine diffeomorphisms of $(X,\omega)$.  Thus,   $\tau \text{PSL}(X, \omega) \tau$ {\em is}  $\text{PSL}(X, \omega)$; compare with Proposition 3.2 of \cite{Mc1}.

In what follows, we will ignore this innocuous conjugation, and consider $\text{PSL}(X, \omega)$ acting on $\mathbb H^2$ by m\"obius transformation.

\subsection{Saddle Connections, Parabolic Directions,  and Special Points}\label{sadParSpecPts} 

Any translation surface of genus greater than one must have singularities.  A geodesic line emanating from a singularity is called a {\em separatrix}, a separatrix connecting singularities (with no singularity on its interior) is a {\em saddle connection}.   We call a point $p \in (X, \omega)$ a {\em connection point} if every separatrix passing through $p$ extends to a saddle connection.   

 We say $\theta$ is a {\em parabolic direction} for  $(X, \omega)$ if there is a parabolic element of $\text{PSL}(X, \omega)$ whose associated affine diffeomorphism preserves the direction $\theta$.   (In fact, the full Veech Dichotomy states that the non-ergodic directions of a Veech surface are parabolic.)    We also call these associated affine diffeomorphisms parabolic.  Expressing $\theta$ as a slope,  the point on $\partial \mathbb H^2$ associated to this slope is $1/\theta$; this is  due to the fact that  $\text{SL}(2, \mathbb R)$ acts on $(X, \omega)$ by way of its {\em linear} action on 
$\mathbb R^2$, but the $\text{SL}(2, \mathbb R)$-action on $\mathbb H^2$ is {\em fractional linear}.   Thus for example, an element of $\text{SL}(2, \mathbb R)$ arising from an element of $\text{Aff}^{+}(X, \omega)$ which fixes the vertical direction (should such an element exist) acts on the extended hyperbolic plane so as to fix the boundary point $\xi =0$, see Figure~1.

A point $p \in (X, \omega)$ of finite orbit under the group of affine diffeomorphisms is called a {\em periodic point}.   A point of finite orbit under the subgroup of affine diffeomorphisms generated by the two (primitive) parabolic diffeomorphisms for  some pair of transverse parabolic directions is called a {\em rational point} for $(X, \omega)$.
In \cite{GHS}, one shows that non-periodic points have dense orbits.   In \cite{HS2}, we show that every connection point is a rational point.

\subsection{Infinitely Generated Groups using Ramified Coverings} 

To date, constructions of infinitely generated Veech groups rely on the following basic fact of 
Fuchsian groups:  (1) a Fuchsian group of the first kind (that is, whose limit set is all of $S^1 =  \partial \mathbb H^2$) is either a lattice group or else is infinitely generated.     In \cite{HS3},   we used two further facts:    (2)  given a  ramified covering $f: X\to Y$ of Riemann surfaces,  upon setting $\omega = f^*\alpha$  the translation surface $(X, \omega)$ is a {\em balanced translation covering}  of the marked translation surface $(Y, \alpha; \text{Br}(f)\,)$ if every point of the fiber above a nonsingular branch point is ramified, in particular,  $\text{PSL}(X, \omega)$ is then commensurable to $\text{PSL}(Y, \alpha; \text{Br}(f)\,)$, where this latter is the  image in  $\text{PSL}(Y, \alpha)$ of the affine diffeomorphisms which stabilize the branch locus of $f$; and, (3)  if there is an element of this branch locus that is non-periodic (that is, such that its $\text{Aff}^{+}(Y, \alpha)$-orbit is of infinite order), then  $\text{PSL}(Y, \alpha; \text{Br}(f)\,)$ is of infinite index in $\text{PSL}(Y, \alpha)$;  in particular, it is thus not a lattice.  Therefore, the commensurable  $\text{PSL}(X, \omega)$ also fails to be a lattice.  

  Given any point of a translation surface, the directions of separatrices through the point are dense (see, say, Lemma 1 of \cite{HS3}, in turn a slight generalization of a result of  \cite{Vo});  the direction of any saddle connection on a Veech surface  is in fact a parabolic direction;  one finds that if 
the branch locus contains a connection point, then the parabolic directions of $\text{PSL}(Y, \alpha; \text{Br}(f)\,)$ are dense in $\partial \mathbb H^2$, and thus the limit set this group is certainly all of $\partial \mathbb H^2$.     Now, if the connection point is non-periodic,  then we conclude that that $\text{PSL}(Y, \alpha; \text{Br}(f)\,)$ is infinitely generated; again, the commensurable $\text{PSL}(X, \omega)$ is therefore also infinitely generated.   
But, we show in \cite{HS3} that there exist $(Y, \alpha)$ that admit non-periodic connection points, and thus proved the existence of infinitely generated $\text{PSL}(X, \omega)$.  \\

 For an explicit example, we showed in \cite{HS3} (see in particular Proposition 7 there) that the genus two surface given by the golden cross can be taken to play the role of $(Y, \alpha)$.  This surface has exactly one singularity, say $s$; let us also choose one of its infinitely many non-periodic connection points, $p$.     If $g=2d$ is even, then we can take $d$ copies of our genus two surface, with a slit from $s$ to $p$; identify these with the usual process after having chosen some transitive permutation on $d$ letters.  This gives a degree $d$ ramified covering, whose branch locus consists of $\{s,p\}$, and which is totally ramified above each of these points.  The Riemann-Hurwitz formula then shows that the covering Riemann surface does have genus $g$.      One can choose the 
 odd genus $g = 3 + 2k$ by way, say, of double coverings:   choose $4k$ connection points of our genus two surface; create $2k$ slits and glue two copies in the standard manner.   For more on creating ramified coverings, see say \cite{D}.

\subsection{Cylinders and Horocyclic Flow}\label{cylHoroFloSect}  

A {\em cylinder} on a translation surface $(X, \omega)$ is a maximal connected set of homotopic simple closed geodesics.  If the genus of $X$ is greater than one, then every cylinder is bounded by saddle connections.   Upon cutting the cylinder along a line perpendicular to the direction of flow, we can represent it as indicated in Figure ~1: of width $w$ and height (or circumference) $h$.   
A cylinder has its affine Dehn twist  
$(x, y) \mapsto (x, \, y + x h/w \bmod h\,)$.   The {\em modulus} of the cylinder is $\mu = w/h$.   If for some direction $\theta$ the surface $(X, \omega)$ decomposes into cylinders in this direction, then we say that $\theta$ is a {\em (completely) periodic direction}.     Veech \cite{Vch1} shows that a completely periodic direction $\theta$  is a parabolic direction if and only if the  various cylinders in the direction $\theta$ gave commensurable inverse moduli:  all $h_i/w_i \cdot w_j/h_j$ are rational.

   The Teichm\"uller {\em horocyclic  flow}  is given by $h_s \circ (X, \omega)$, where
 $h_s := \begin{pmatrix} 1&0\\ s&1\end{pmatrix}$, with $s\in \mathbb R$. The orbit of  $(X, \omega)$ under this flow is called a {\em horocycle}.

Suppose that $(X, \omega)$  is completely periodic in the vertical direction, that is 
$(X, \omega)$ has a decomposition into vertical cylinders.  The effect of  $h_s$, with $s>0$, restricted to  a cylinder of the decomposition is to change the relative  positions of the singularities on the right hand side of the cylinder with respect to those on the left.   It is easily seen that on a given vertical cylinder of height $h$ and width $w$, that a standard affine twist on the cylinder has been performed when $s = h/w$.   (Thus, when $s$ is the multiplicative inverse of the modulus of the cylinder.)   See Figure~1.

\bigskip 

\centerline{
\vbox{
\includegraphics{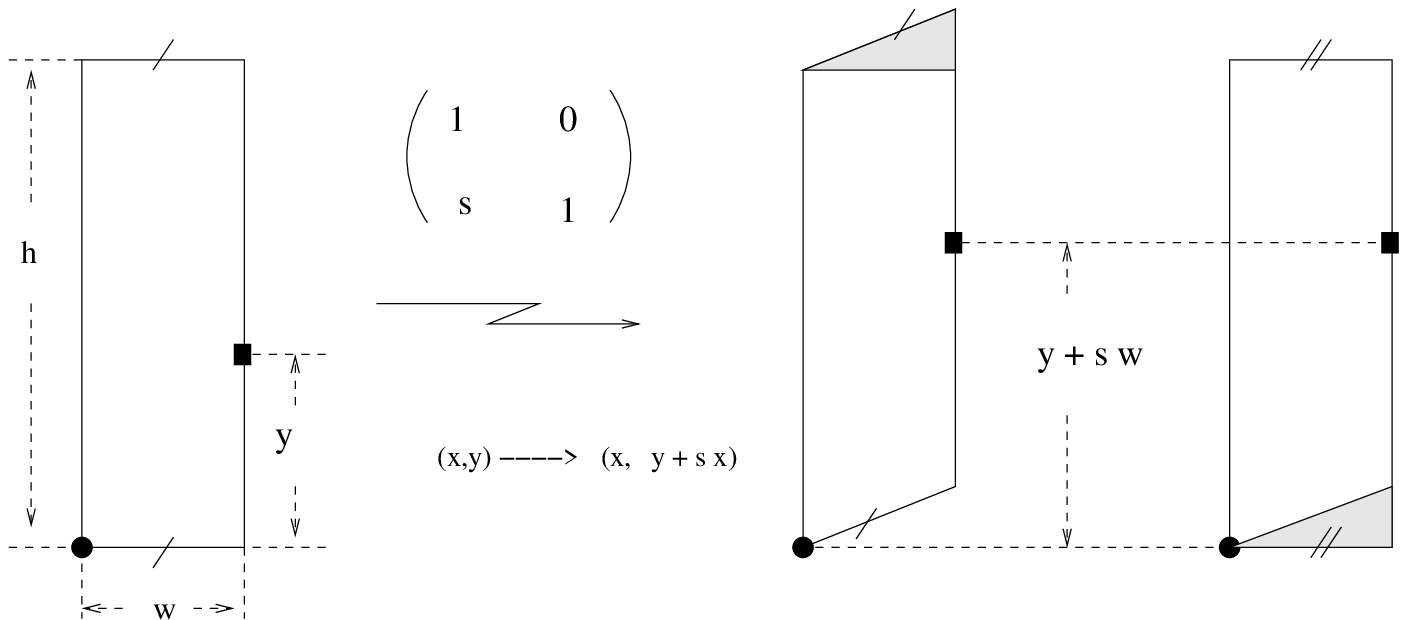}
}}
\vspace{ .3 cm}
\centerline{
{Figure 1.}  Horocyclic flow:   moving a singularity.}
\bigskip

If the various moduli of the cylinders are commensurable, then there is a least positive $s$ such that $h_s$  induces on each cylinder an integer power of the basic 
Dehn twist of the cylinder.      This 
$h_s$ is then easily seen to belong to $\text{SL}(X, \omega)$.     Indeed, the element $h_s$ then generates a parabolic subgroup of  $\text{PSL}(X, \omega)$, and the   $\text{PSL}(X, \omega)$-conjugacy class of this parabolic subgroup corresponds to a {\em cusp}.

\subsection{Geodesic Flow, Horodisks and Fat Directions} \label{geoFloHoroFat} 

Teichm\"uller  {\em geodesic flow} is given by $g_t \circ (X, \omega)$, where $g_t = \begin{pmatrix} e^{t/2}&0\\0&e^{-t/2}\end{pmatrix}$, with $t\in \mathbb R$.   For positive $t$, the effect of $g_t$ restricted to a vertical cylinder decomposition is to decrease heights and increase widths.    Note that  vertical saddle connection lengths tend to zero as $t \to \infty$.  Thus, if $ (X, \omega)$ has saddle connections in the vertical direction, then for  $t$ sufficiently large,  the shortest saddle connections for $g_t \circ (X, \omega)$ are in the vertical direction.        

A {\em horodisk} of index $T$ corresponding to a direction $\theta$ of $(X, \omega)$ is the union over $t>T$ of all the horocycles $h_s (\,g_t \circ (X, \omega)\,)\,$.      

\begin{defin}  A direction $\theta$ of  $(X, \omega)$  is a {\em fat direction}  if there is a sufficiently large positive real $T$ such no two points of the  horodisk of index $T$ corresponding to $\theta$ are  $\text{SL}(X, \omega)$-equivalent.   
\end{defin} 

\subsection{Ends of Hyperbolic Surfaces}\label{endHypSurf} 
An {\em end} of a hyperbolic manifold $\mathcal S$ is determined by a cover of $\mathcal S$ by compact sets $K_i$, with each $K_i$ contained in the interior of $K_{i+1}$; an end $\mathcal E$  is a sequence of components $\mathcal E_i$ in the complement of $K_i$ such that $\mathcal E_{i+1} \subset \mathcal E_i$.   The definition of an end is independent of choice of sequence $K_i$, see say \cite{ThurstonBook}.  
An end $\mathcal E$  of a hyperbolic surface is called a {\em cusp} if there is a subset $\mathcal P \subset \mathcal S$ that is conformally equivalent to the complex punctured unit disk $\{z\,|\, 0 < |z|<1\,\}$ such that for all sufficiently large $i$, $\mathcal E_i\subset \mathcal P$.   An end is called a hole if there is a subset $\mathcal A \subset \mathcal S$ that is conformally equivalent to a complex annulus  $\{z\,|\, 1 < |z|<r\,\}$ with $r>1$ such that for all sufficiently large $i$, $\mathcal E_i\subset \mathcal A$.   These two types of ends are called finite, an {\em infinite end} is an end which is not finite.    In a terminology introduced by Haas \cite{Ha}, one differentiates between ends of the first kind and infinite ends of the second kind; the infinite ends that we discuss are of the first kind.
To be precise, we say that    $\mathcal E$ is an {\em infinite end of the first kind} if 
there is a subset $\mathcal U \subset \mathcal S$ that is conformally equivalent to the complex upper half-plane   $\{z\,|\, \text{Im}(z)>0\,\}$  such that for all sufficiently large $i$, $\mathcal E_i\subset \mathcal U$. Note that each of the above definitions  continues to be valid upon allowing $\mathcal S$ to be a hyperbolic orbifold, although there is then also the further possibility that the set of elliptic fixed points enters arbitrarily far into an end.\\

\begin{defin}  Given a Fuchsian group $\Gamma$, and a point $\xi \in \partial \mathbb H^2$, we say that $\xi$ {\em determines an infinite end} of $\Gamma \backslash \mathbb H^2$ if there is some horodisk based at $\xi$ such that any $\Gamma$-orbit in $\mathbb H^2$ meets this horodisk in at most one point.    Equivalently,  there is an infinite end of 
$\Gamma \backslash \mathbb H^2$ which lifts to include a horodisk based at $\xi$.   Analogously, we say that  a parabolic fixed point $\xi$ of $\Gamma$  {\em determines a cusp}. 
\end{defin} 
 
  \bigskip 

 The following two lemmas are virtually immediate.   
 
 \begin{lem} \label{niceEnds}     If the Fuchsian group $\Gamma$ is of the first kind, then any end of  $\Gamma \backslash \mathbb H^2$ is either a cusp or an infinite end as above.   
 \end{lem}

 \begin{lem} \label{fatPtsAreInfEndDir}    Let $(X, \omega)\in \Omega \mathcal M_g$, and let   $\theta\in \mathbb P^1(\mathbb R)$ be a direction for $(X, \omega)$.   Let $\xi \in \partial \mathbb H^2$ be given by $\xi = 1/\theta$.  If $\theta$ is a parabolic direction, then $\xi$ determines a cusp of $\text{PSL}(X, \omega) \backslash \mathbb H^2$.  If $\theta$ is a fat direction, then  $\xi$  determines an infinite end of $\text{PSL}(X, \omega) \backslash \mathbb H^2$.   \end{lem}

 \bigskip 
 For ease of discussion, we introduce the following terms.  

\begin{defin}  Given a (possibly ramified) covering of hyperbolic surfaces or orbifolds $\pi: \mathcal S \to \mathcal T$, we say $\pi$ is {\em balanced with respect to ends} if the 
image of any end of $\mathcal S$ is an end of $\mathcal T$ and also the pre-image of any end of $\mathcal T$ is an end of $\mathcal S$.   The terms  {\em balanced with respect to infinite ends} and {\em balanced with respect to cusps} are analogously defined.
\end{defin} 
 
 \bigskip 
 
The following is well-known.

\begin{lem}  \label{finCoveringEnds}  Finite degree ramified coverings of hyperbolic surfaces or  orbifolds   are balanced with respect to each of the following: 
ends, infinite ends, and cusps.
\end{lem}  

 \bigskip 
Proposition 3 of \cite{HS3}  can be reworded as the following.   

\begin{prop}\label{ourEnds} Let $(Y, \alpha)$ be a Veech translation surface. 
Suppose that  $p\in Y$ is a non-periodic connection point.   Let $\text{PSL}(Y, \omega; p)$ denote the subgroup of $\text{PSL}(X, \omega)$ arising from the stabilizer of $p$ in $\text{Aff}^{+}(X, \omega)$.   Then the  (possibly ramified) covering $\pi: \text{PSL}(Y, \omega; p)\backslash \mathbb H^2 \to \text{PSL}(Y, \omega)\backslash \mathbb H^2$  is balanced with respect to ends.     However, this covering is not balanced with respect to either cusps or infinite ends; more precisely:  $\text{PSL}(Y, \omega; p)\backslash \mathbb H^2$ has infinite ends, but $\text{PSL}(Y, \alpha)\backslash \mathbb H^2$ has no infinite ends. 
\end{prop}  

\begin{rem} Points defining our ends are Dirichlet points.  Recall that a {\em Dirichlet point} is a point $\xi \in \mathbb H^2 \cup \partial \mathbb H^2$ which is on the boundary of some Dirichlet region --- thus, $\xi$ admits some $a \in \mathbb H^2$ such that $\xi$ is on the boundary of the Dirichlet fundamental region for $\Gamma$.   
One can show, see \cite{N1, N2}, that a Dirichlet point $\xi \in \partial \mathbb H^2$ is 
either a parabolic fixed point of $\Gamma$, or is such that the $\Gamma$-orbit of $a$ meets the boundary of a horodisk based at $\xi$ whose interior contains no element of this orbit.   

Recall that $\xi \in \partial \mathbb H^2$ is called, after Sullivan \cite{S}, a {\em Garnett point} 
for $\Gamma$ relative to $a$ if the largest  horodisk  based at $\xi$ whose interior contains no element of  the $\Gamma$-orbit of $a$, also does not meet this orbit.   A point $\xi$ is a {\em rigid Garnett point} for $\Gamma$ if $\xi$ is Garnett with respect to all $a\in \mathbb H^2$.
Now, if  $G$ is any  Fuchsian group, and   $\xi$ is a parabolic point for $G$, then for any subgroup $\Gamma \subset G$ admitting a fundamental domain which contains some non-trivial horodisk based at $\xi$, one can easily show that there is some $a \in \mathbb H^2$ such that $\xi$ is a Dirichlet point for $\Gamma$ relative to $a$.   Thus,  the points which determine our infinite ends are not  rigid Garnett points.   We do not yet know if they can in fact be Garnett points. 
\end{rem}

\section{Projecting Horocycles to $\mathcal M_g$}\label{projHoro} 

We discuss the projection of certain horocycles of the Teichm\"uller disk of some $(X, \omega)$ to $\mathcal M_g$.  

\subsection{Factoring Horocycles Through Tori} 
We discuss the avatars  in the setting of $\Omega \mathcal M_g$ of certain elements of the Teichm\"uller modular group.

Suppose  that $(X, \omega)$ admits a cylinder $\mathcal C$ in the vertical direction.     Consider  the basic affine Dehn
 twist on  $\mathcal C$,  expressed  in standard coordinates simply as $(x, y) \mapsto (x, \, y + x h/w \bmod h\,)$, where $h,w$ are respectively the height and width of $\mathcal C$, see Figure ~1.    Recall that $\mu := w/h$ is defined as the {\em modulus} of the cylinder.  This Dehn twist  acts as the identity on the vertical boundaries of the cylinder; we can extend by the identity map to obtain a well-defined homeomorphism, $\sigma$, of the underlying topological surface.  But, $\sigma$ allows us to give a differentiable structure on the receiving surface (such that the image of  the singularities are again singularities); furthermore, $\sigma$ has locally constant derivative.   Hence,   using $\sigma^{-1}$ we pull-back the translation atlas from $(X, \omega)$ to give a new translation structure. Now,  each of the various cylinders $\mathcal C_i$ has image under $\sigma$ of equal modulus;  by standard results in complex analysis, see say \cite{G}, $\sigma$ is conformal.   Therefore, the complex structure pulled-back by $\sigma^{-1}$ is the same as that of $X$.   
Thus, we find a mapping of translation surfaces,  $\sigma: (X, \omega) \to (X, \beta)\,$.
(Of course, $\omega=\beta$ only when  $\mathcal C$  is all of $X$.)  
We thus have that  $\pi \circ \sigma (\,  (X, \omega)\,) = \pi (\,  (X, \omega)\,)$, where 
$\pi: \Omega \mathcal M_g \to \mathcal M_g$ is the projection.  Note that the horocyclic flow $\{\, h_s \circ (X, \omega)\,|\, s\in \mathbb R\,\}$ restricts to induce a map on  $\mathcal C$ which commutes with $\sigma$ acting on $\mathcal C\,$;  modulo the action of  $\sigma$,  this map is  parametrized by $s \bmod  \mu^{-1}\, $.

  Now suppose that there are $n$ cylinders in the vertical direction, $\mathcal C_1$, $\dots$, $\mathcal C_n$.  Let $\mu_i$  denote the modulus, and let $\sigma_i$ denote the  map as above,  corresponding to the cylinder $\mathcal C_i\,$;  the various $\sigma_i$  commute, as maps on the underlying surface $X$.    Similarly, we can define an $\mathbb R^n$-action on the underlying surface: given $s = (s_1, \dots, s_n) \in \mathbb R^n$, in cylinder $\mathcal C_i$, we twist by amount $s_i$; outside the collection of cylinders, the action is defined to be trivial.    Since this action commutes with the action of the $\sigma_i$, we find that the projection to $\mathcal M_g$ of 
 $\mathbb R^n\,\circ (X, \omega)$ factors through  $\mathcal T_{X,\omega} := 
 \mathbb R^n\,\circ (X, \omega)\}/\langle \sigma_i\,\rangle\,.$  But, $\mathcal T_{X,\omega}$ is isomorphic to the $n$-torus, $\oplus_{i=1}^n\;  \mathbb R/\mu_{i}^{-1} \mathbb Z\,$.
 
   If in fact $(X,\omega)$ admits a cylinder decomposition given by the various $\mathcal C_i$,  then the horocycle  $\{\, h_s \circ (X, \omega)\,|\, s\in \mathbb R\,\}$ lies in  $\mathbb R^n\,\circ (X, \omega)$.   Indeed, the projection of this horocycle to $\mathcal M_g$ factors through a curve on $\mathcal T_{X,\omega}$.     To generalize slightly, if   $\theta$  denotes the direction of a cylinder decomposition of $(X, \omega)$ into $n$ cylinders,  we have that the projection of the corresponding horocycle factors through a curve on the $n$-torus denoted $\mathcal T_{X,\omega, \theta}$. 

Commensurability of non-zero real numbers is an equivalence relation; for each commensurability class of the moduli $\mu_i$, there is a corresponding smallest positive value of $s$ which induces on each of the corresponding cylinders $\mathcal C_i$ an integer power of its basic Dehn twist.       Let $m$ be the number of these classes; for the $j$-th class, let $\tilde \mu_{j}^{-1}$ be the least common positive integer   multiple of the various $\mu_{i}^{-1}$ in this class.   We then let $\tilde \sigma_j$ be the composition of the various $\sigma_{i}^{r(i)}$ in this class, where $r(i)$ is such that $r(i) \cdot \mu_{i}^{-1}=\tilde \mu_{j}^{-1}$.  Again, the projection of the horocycle factors through the torus corresponding to the  group generated by the $\tilde \sigma_j$.   We find that in fact our horocycle factors through a curve on the $m$-torus $\mathcal{\tilde T}_{X,\omega, \theta} \cong \oplus_{j=1}^m\;  \mathbb R/\tilde \mu_{j}^{-1} \mathbb Z\,$.    

Of course,  $\theta$ is a parabolic direction if and only if there is a single commensurability class, in which case: $\tilde \sigma_{1} \in \text{Aff}^{+}(X, \omega)$, and  the horocycle factors through a circle; this circle lies on (a finite component  of the sequence which defines the end which is)   the corresponding cusp of $\text{PSL}(X,\omega)\backslash \mathbb H^2$.    In all other cases,  there is more than one commensurability class defining $\mathcal{\tilde T}_{X,\omega, \theta}$,    and  just as $\{ (s \bmod \tilde\mu_{1}^{-1}, \dots, s \bmod \tilde\mu_{m}^{-1})\,|\, s\in \mathbb R\,\}$ is dense in $\oplus_{j=1}^m\;  \mathbb R/\tilde \mu_{j}^{-1} \mathbb Z\,$, 
so is the image of the horocycle is dense on $\mathcal{\tilde T}_{X,\omega, \theta}\,$.

\subsection{Geodesic Flow and Horocyclic Projections} \label{geoFloHoroPro} 

 For positive $t$, the effect of 
$g_t = \begin{pmatrix} e^{t/2}&0\\0&e^{-t/2}\end{pmatrix}$ on a vertical cylinder decomposition of $(X, \omega)$ is to preserve the number of cylinders, while decreasing heights and increasing widths.  Indeed, $g_t$ acts so as to multiply vertical cylinder moduli by $e^{t}$.   Therefore,  if $\theta$ is a direction of a cylinder decomposition of $(X, \omega)$ as above, then $\mathcal{\tilde T}_{g_t\circ(X,\omega), \theta}$     has the same dimension as $\mathcal{\tilde T}_{(X,\omega), \theta}\,$.

 Note that  vertical saddle connection lengths tend to zero as $t \to \infty$; indeed (having begun with vertical saddle connections),  for $t$ sufficiently large,  the shortest saddle connections for $g_t \circ (X, \omega)$ are in the vertical direction.   But then,  $g_t \circ (X, \omega)$ leaves every compact  subset of $\mathcal M_g$.    Hence, as $t \to \infty$, $g_t \circ (X, \omega)$ penetrates into an end of $\mathcal M_g$, see say \cite{MT}.

\section{Some Elementary Tools} 

In \cite{HS2}, we studied an invariant of the geometry at a cusp of a Teichm\"uller complex geodesic --- the multiplicities and relative lengths of saddle connection vectors.   We use the following elementary notion in a similar manner. 

\begin{defin} If $p$ lies interior to the sides of a cylinder, then the {\em splitting ratio} of $p$ is the ratio 
of the distance from  the left hand side of the cylinder to  $p$ to the width of the cylinder.   \end{defin} 

Lemma 4 of \cite{HS1} can be reworded as the following.

 \begin{lem} \label{irrSplit}   Let $(Y, \alpha)$ be a translation surface.  If $\theta$ is a completely periodic direction and $p\in Y$ has irrational splitting ratio for some cylinder of the cylinder decomposition in the direction $\theta$, then $\theta$ is a fat direction for 
  $(Y, \alpha; p)$.   
  \end{lem}  

 \begin{lem} \label{splitRatioFixed}    Splitting ratios are $\text{SL}(2, \mathbb R)$ invariant.
 \end{lem}  

\noindent 
{\bf Proof}    Let $(X, \omega)$ be a translation surface, and suppose that $p\in X$ splits some cylinder with ratio $\theta$.   For ease of discussion, we may and do assume that the cylinder, say $\mathcal C$, is a vertical cylinder.   Suppose that  $A \in \text{SL}(2, \mathbb R)$.   

We develop   $\mathcal C$ to $\mathbb R^2$ in the standard fashion.   
 The ratio $\theta$ is the ratio of the length of the horizontal line segment passing from 
 the left hand side of $\mathcal C$ to $p$ to the width of the cylinder.  Of course, this second length is the length of the horizontal line segment joining the sides of $\mathcal C$ and passing through $p$.

The image of the cylinder by $A$ is again a cylinder, say $\mathcal D$, with obvious development.  
The horizontal line segment crossing $\mathcal C$ upon which $p$ lies has as its image some segment in $\mathcal D$;  since $A$ acts linearly,  the ratio of the length of the initial segment of this image --- from the appropriate side of the cylinder to the image of $p$ --- to the length the segment  is again $\theta$.

Draw the line segment perpendicular to the sides of the cylinder which passes through the image of $p$.  Due to the equality of the opposite angles so defined,  ratios of corresponding sides of the triangles so formed are equal.   A trivial cancellation, see Figure 2, shows that the splitting ratio is indeed preserved. 
\qed

\bigskip 

\centerline{
\vbox{
\includegraphics{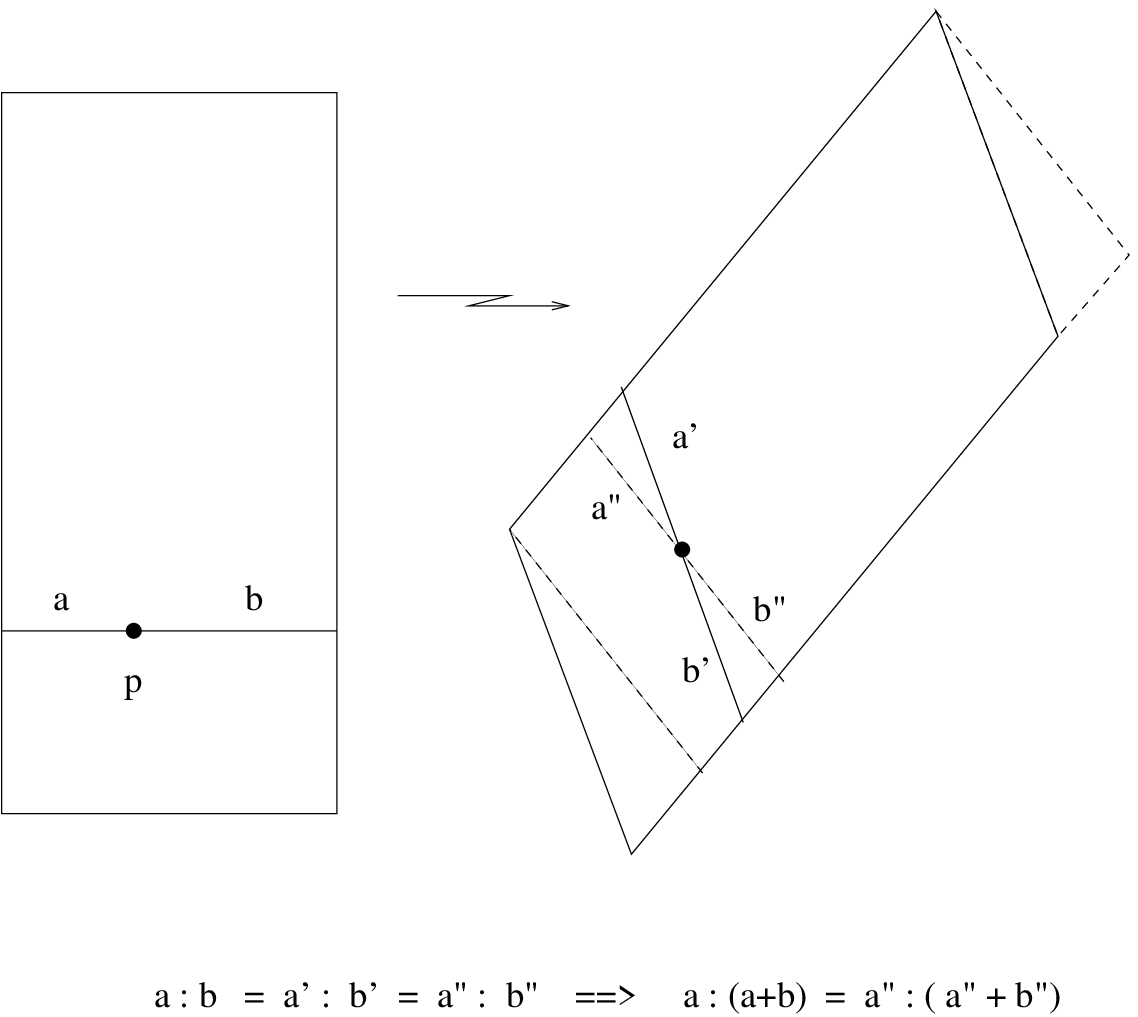}
}}
\vspace{ .3 cm}
\centerline{
{Figure 2.}  Preserving splitting ratios.}

\bigskip 

We have the immediate consequence. 

 \begin{corol} \label{affPresSplitt}   Affine diffeomorphisms preserve splitting ratios of
 cylinders.   
\end{corol}  

We distinguish infinite ends by use of  the following elementary well-known result in the theory of continued fractions.   We use a flat notation for continued fractions; thus, for example $[a; b, c+ \alpha ] = a+ \dfrac{1}{b+ \dfrac{1}{c + \alpha}}\,$.     We label the partial quotients of a real number $\theta \in [0,1]$ by $0, a_1$, $a_2$, $\dots$.      These give rise to the convergents to $\theta$:   $p_n/q_n = [0;a_1, a_2, \dots, a_n]$.    Let $\{x\}$ denote the fractional part of a real number: $\{x\} = x - \lfloor \,x\,\rfloor\,$.  

 \begin{lem} \label{rotateIrrationals}   Suppose that $\mathcal C$ and $\mathcal D$ are cylinders lying in tranverse directions which have non-trivial intersection.  If $p$ irrationally splits $\mathcal C$, then the orbit of $p$ under the affine Dehn twist of $\mathcal C$ meets the cylinder $\mathcal D$ so as to give infinitely many distinct irrational splitting ratios.   
\end{lem}  

\noindent 
{\bf Proof }   Let  $\theta$ denote the splitting ratio of $p$ with respect to $\mathcal C$.   By Lemma ~\ref{splitRatioFixed}, we  may normalize such that $\mathcal C$ is in the vertical direction and $\mathcal D$ is in the horizontal direction.    Let $\mathcal R$ be a rectangle in which our two cylinders meet.     Let us develop  the cylinders such that $p$ lies on the $x$-axis, and such there are $y_0$ and $\lambda$ such that $\mathcal R$ is bounded below by the line $y= y_0$ and above by $y= y_0 + \lambda$.    Let us denote by $h$ the vertical height $\mathcal C$.    Thus in these coordinates, the orbit of $p$ by the vertical Dehn twist  has heights $n \theta h  \, \bmod h\,,$ with $n \in \mathbb N$.     But, this is exactly the set of $\{n \theta\} h\,$.

It suffices for us to show that
there are infinitely many distinct ratios  $(\, \{n \theta\} h\,)/\lambda$ which are irrational, and lie between $y_0$ and $y_0 + \lambda$.     For completeness' sake, we show this elementary, presumably well-known, result.

 Let $\nu := \lambda/h$.      Now, if 
 $\{n \theta\}$ is a rational multiple of $\nu$, then it is necessary  that $\theta = \dfrac{a \nu + b}{c}$ for $a,b,c$ integers, with $a$, $c$ nonzero.     
 
Suppose now that   $\theta$ is of the form   $\dfrac{a \nu + b}{c}$,  and there is an $n$ such that $\{n \theta\} = r \nu/s$ for nonzero integers $r$ and $s$.     We find:  
\[n \dfrac{a \nu + b}{c} - k = r \nu/s\,,\; \text{with}\,  k = \lfloor\, n \theta\,\rfloor\,.\]
  Since $\theta$ is irrational, it is not linearly related to a rational number.  Hence,    in particular, we must have that $k = nb/c$.   That is,   $\lfloor n \theta \, \rfloor = n\, b/c$; but  then  $0<\theta - b/c<1/n\,$.  This can only occur for a finite number of $n$.

  Since in the worst case scenario there are only a finite number of orbit elements to avoid, and the orbit of $p$ is dense in the circle on the vertical cylinder, we conclude that there are always infinitely many elements of the orbit which irrationally split the horizontal cylinder. \qed

\begin{lem}  \label{connPtsForever}  The affine diffeomorphisms of a translation surface preserve the set of connection points.  \end{lem} 

\noindent 
{\bf Proof}    Suppose that $\phi \in \text{ Aff}^{+}(X, \omega)$ is an affine diffeomorphism and $p\in X$ is a connection point.   Given a separatrix passing from a singularity to $\phi(p)$, we apply $\phi^{-1}$ to find a separatrix passing from some singularity through $p$.   Since $p$ is a connection point, this separatrix extends to a saddle connection.   But, the image by $\phi$ of this extended separatrix extends our initial separatrix passing through $\phi(p)$.  Since $\phi$ sends saddle connections to saddle connections, we are done. \qed 

\bigskip 

\section{Infinitely Many Cusps}\label{manyCusps} 
We use the invariant studied in \cite{HS2}  to distinguish  cusps.    

\begin{prop}  \label{infCusps}  Let $(X, \omega)$ be a translation surface which is a finite ramified covering of a Veech surface $(Y, \alpha)$, with ramification above some non-periodic connection point and possibly the singularities of $(Y, \alpha)$.    Then  the Teichm\"uller complex geodesic  $\mathcal G(X, \omega)$ has infinitely many cusps.   
\end{prop}  

\noindent 
{\bf Proof}   The translation structure of $(X, \omega)$ arises as the pull-back by a ramified holomorphic map $f: X \to Y$ of Riemann surfaces, with $\omega = f^{*}\alpha$.     Let $p$ denote the unique non-periodic connection point in the branch locus.  The Veech group $\text{PSL}(X, \omega)$  is commensurable to the subgroup of $\text{PSL}(Y, \alpha)$ formed by the derivatives of the  stabilizer in $\text{Aff}^{+}(X,\alpha)$ of $p$.

By Theorem 2 of \cite{GHS}, the image of the non-periodic $p$ under the action of the affine diffeomorphism group is dense on $(Y,\alpha)$.   Let $l= l_{\text{min}}(Y,\alpha)$ be the length of the shortest saddle connection of  $(Y,\alpha)$.     Given $\epsilon < l/2$,      we can find an image of $p$  under $\text{Aff}^+(X, \omega)$, say $q$,  which is within $\epsilon$ of one of the singularities, say $\sigma$.      There is a unique 
separatrix of length less than $\epsilon$ connecting $q$ and $\sigma$.  
 By Lemma ~\ref{connPtsForever}, we can extend  this separatrix to a saddle connection, passing from $\sigma$ to a singularity $\tau$. (Of course, $\tau$ may equal $\sigma$.) 

Let $\phi$ be an affine diffeomorphism mapping $q$ to $p$.  The image by $\phi$ of the above saddle connection is a saddle connection passing through $p$.   Since $(Y,\alpha)$ is a Veech surface, the direction of any saddle connection is a parabolic direction.   Since $p$ lies on this saddle connection,  $p$ is preserved by  the Dehn twist of either of the cylinders (in this given direction) upon which $p$ lies.     Thus this direction is a parabolic direction for $(X,\omega)$.     

Consider a single  lift to $(X,\omega)$ of the saddle connection that we have just determined.   This lift is the union of two saddle connections on $(X,\omega)$:  the first joins a lift of $\phi(\sigma)$ to a lift of $p$,   the second lies on another branch of the (ramified) covering space  $X$, joining the lift of $p$ to a lift of $\phi(\tau)$.    (The connection vectors of these saddle connections are parallel!)

Now, the lengths of these  saddle connections on  $(X,\omega)$ are equal to the lengths 
of the corresponding segments on  $(Y,\alpha)$.   But,  the affine diffeomorphism $\phi$  preserves {\em ratios} of segment lengths;  thus the ratio of the lengths of the saddle connections on  $(X,\omega)$ equals the ratio of the length along the separatrix from $\sigma$ to $q$ to the length along its continuation from $q$ to $\tau$.   We find that this ratio is less than $\epsilon/ 2 l$. 

Our construction thus allows us to find pairs of parallel connection vectors, in parabolic directions, of $(X,\omega)$ that are of arbitrarily small ratio.   However, the number of connection vectors in any given direction is finite.   Furthermore, the ratios of the various connection vectors lying in a parabolic direction is invariant under the action of $\text{SL}(X, \omega)$ --- see \cite{HS2}.   We conclude that $(X, \omega)$ has infinitely many inequivalent parabolic directions.   But, each of these parabolic directions determines a distinct cusp of 
$\text{PSL}(X, \omega)\backslash \mathbb H^2$.   Since the Teichm\"uller complex geodesic  $\mathcal G(X, \omega)$ is a generically injective image of this quotient, it too has infinitely many cusps.
\qed 

\bigskip 

\section{Infinitely Many Infinite Ends}\label{manyEnds}

We now give the proof of Theorem \ref{infFatBoys}. 
Thus, we continue with $(X,\omega)$ a ramified covering  of $(Y,\alpha)$ as above. 
The periodic directions on $(X,\omega)$ are exactly those given by the periodic, and hence parabolic, directions of $(Y,\alpha)$.      We show that $(X,\omega)$ has many directions which are periodic, but not parabolic.

\begin{prop}  \label{infEnds}  Let $(X, \omega)$ be a translation surface which is a finite ramified covering of a Veech surface $(Y, \alpha)$, with ramification above some non-periodic connection point and possibly the singularities of $(Y, \alpha)$.    Then  the Teichm\"uller complex geodesic $\mathcal G(X, \omega)$ has infinitely many infinite ends.   
\end{prop}  

\noindent 
{\bf Proof}   Let $p$ be the branch point, thus $p$ is a non-periodic connection point of $(Y, \alpha)$. 
There is thus some parabolic direction of $(Y,\alpha)$ for which $p$ irrationally splits a cylinder.     Without loss of generality, we may normalize so that $p$ irrationally splits a vertical cylinder and such that the horizontal is a parabolic direction for     $(Y, \alpha)$.      

Consider  the vertical cylinder $\mathcal C$ in which $p$ lies, and any horizontal cylinder $\mathcal D$ which meets $\mathcal C$.     By Lemma ~\ref{rotateIrrationals},  the set of splitting ratios of $\mathcal D$ determined by the orbit of $p$ under the vertical affine Dehn twist includes an infinite set of irrational 
numbers.      There is some finite power of this Dehn twist which is induced by an element $\phi$ of the affine diffeomorphism group of $(Y, \alpha)$; the orbit of $p$ under the subgroup generated by this multiply-composed Dehn twist still provides  an infinite set of distinct irrational splitting ratios.

For each image $q= \phi^n(p)$ of this latter type, apply the corresponding  inverse $\phi^{-n}$.   Now,
$\phi^{-n}(q)=p\,$; $\,\phi^{-n}(\mathcal D)$ is a cylinder in a corresponding parabolic direction for $(Y, \alpha)\,$;  and, $p$ splits $\phi^{-n}(\mathcal D)$ with the ratio that $q$ splits $\mathcal D$.    Of course, given any direction, $p$ lies interior to at most one cylinder of the possible decomposition in that direction;   thus the splitting ratio defined by $p$ in that direction is well-defined.       Furthermore, the splitting ratio that $p$ gives is constant on directions equivalent under the stabilizer of $p$ amongst the affine diffeomorphisms of $(Y, \alpha)$.    Since the $\text{Aff}^{+}(Y, \alpha)$-orbit of $p$ gives infinitely many  distinct irrational splitting ratios, we conclude that there are infinitely many  
 fat directions that are inequivalent  under the stabilizer of $p$.   
 
\bigskip 

Therefore,  $\text{PSL}(Y, \alpha; p)\backslash \mathbb H^2$ has 
 infinitely many infinite ends.   But,  $\text{PSL}(X, \omega)$ is commensurable with 
$\text{PSL}(Y, \alpha; p)$; hence, $\text{PSL}(Y, \alpha; p)\backslash \mathbb H^2$  and $\text{PSL}(X, \omega)\backslash \mathbb H^2$ admit a common finite degree (possibly ramified) covering.   Since such coverings are balanced with respect to infinite ends, $\text{PSL}(X, \omega)\backslash \mathbb H^2$ also has infinitely many infinite ends. 
But, the Teichm\"uller complex geodesic $\mathcal G(X, \omega)$ is a generically injective image of $\text{PSL}(X, \omega)\backslash \mathbb H^2$, it too must have infinitely many infinite ends. 
\qed   

\bigskip 

We next observe that the above proof actually shows a bit more. 
\begin{lem}  \label{fatLimits}  With notation as above, let $\theta$ be a fat direction of 
$(X, \omega)$.   Then $\theta$ is the limit of $\text{PSL}(X, \omega)$-inequivalent fat directions of $(X, \omega)$.
\end{lem}  

\noindent 
{\bf Proof}  From construction, $\theta$ is a parabolic direction for $(Y, \alpha)$ and 
the non-periodic connection branch point $p$ irrationally splits a cylinder $\mathcal C$ of $(Y, \alpha)$ in this direction.  Recall that some power of the Dehn twist in this cylinder is induced by an element $\phi$ of the affine diffeomorphism group of $(Y, \alpha)$.   
Choose any parabolic direction for $(Y, \alpha)$ that is transverse  to $\theta$.   There is some cylinder $\mathcal D$ in this direction that meets $\mathcal C$.   The previous proof shows that there is a subsequence of $\phi^{-n}(\mathcal D)$ giving cylinders which are split by $p$ with distinct irrational splitting ratios.    Since $\phi$ is parabolic, preserving the direction $\theta$, the direction of the cylinder $\phi^{-n}(\mathcal D)$ tends to $\theta$ as $n$ tends to infinity.   Thus, $\theta$ is the limit of these inequivalent fat directions for $(Y, \alpha; p)$.  Finally, commensurability of groups allows us to choose a subsequence of these directions consisting of inequivalent fat directions for $(X, \omega)$.
\qed

\section{Locating Ends}\label{locEnd}
 
 We now suppose $f: (X,\omega) \to (Y, \alpha)$ is such that $(X, \omega)$ satisfies the hypotheses of Theorem \ref{infFatBoys}; thus: $(Y, \alpha)$ is 
a Veech surface $(Y, \alpha)$ with distinguished non-periodic connection point $p$,  $\omega = f^{*} \alpha$, and the branch locus of $f$ is  $p$ and possibly singularities of $(Y, \alpha)$.     The case of ramification above multiple non-periodic connection points is similar; for clarity of exposition, we content ourselves with this main case.

  By construction, $\text{SL}(X, \omega)$ is commensurable to   $\text{SL}(Y, \alpha; p)$.  
 For clarity of discussion, {\em suppose further} that the marked translation surface $(Y, \alpha;p)$ has some connection vector which is uniquely realized.   This is easily arranged by choice of $p$.   With this added hypothesis,   Proposition 5.6 of \cite{Vo}  applies to show that $\text{SL}(X, \omega) \subset  \text{SL}(Y, \alpha; p)$.

\subsection{Orbits of Coverings}\label{orbCov}
 
 Our basic tool is the following.

\begin{lem}  \label{everyAffHasBldCover}  Let $f: (X, \omega)\to (Y, \alpha; p)$ be a balanced covering of translation surfaces as above.  For each $\phi\in \text{Aff}^{+}(\,(Y, \alpha)\,)$, there is some $(Z,\beta)$ with $\pi(\,(Z,\beta)\,) \in \mathcal G(X, \omega)$  such that $\phi\circ f: (Z, \beta)\to (Y, \alpha; \phi(p)\,)$ is a balanced translation covering. 
\end{lem}  

\noindent 
{\bf Proof}  We take $A\in \text{SL}(Y,\alpha)$ such that $\phi$ gives an isomorphism of translation surfaces of $A\circ(Y,\alpha)$ with $(Y, \alpha)$.    Let us now excise the singularities from $ A\circ(X, \omega)$ to form $ A\circ(X, \omega)'$, and excise both $\phi(p)$ and the singularities from $(Y,\alpha)$ to form $(Y,\alpha)''$.  We then find that  $\phi\circ f: A\circ(X, \omega)' \to (Y,\alpha)''$ is a covering of topological surfaces.   Letting $(Y,\alpha)''$ have its induced translation structure from  $(Y,\alpha)$, we can pull this back by  $\phi\circ f$ to 
induce a possibly new structure on $A\circ(X, \omega)'$;  we can complete this flat metric to form some $(Z, \beta)$.  By construction,   $\phi\circ f$ gives a balanced translation covering by $(Z,\beta)$.   The translation structure of $(Z, \beta)$ is compatible with a unique complex structure; this complex structure is also (after possibly extending  over removable singularities) the completion of the pull-back of the complex structure on    $(Y,\alpha)''$,  which in turn is simply the restriction of the complex structure on $(Y,\alpha)$.   But, $\phi$ does not change complex structure, hence the pulled-back complex structure is exactly that of $A\circ(X, \omega)$.  Therefore,   $\pi(\, (Z, \beta)\,) \in \mathcal G(X, \omega)$.   
\qed

\bigskip

\begin{lem}  \label{lineUpCovs}  Let $f: (X, \omega)\to (Y, \alpha; p)$ be a balanced covering of translation surfaces as above.  Then for all  $A \in \text{SL}(Y, \alpha)$, there is some $(Z,\beta)$ with $\pi(\,(Z,\beta)\,) = \pi(\,A\circ(X, \omega)\,)$ and some $\phi\in \text{Aff}^{+}(Y, \alpha)$ such that $\phi\circ f: (Z, \beta)\to (Y, \alpha; \phi(p)\,)$ is a balanced translation covering.    
\end{lem}  

\noindent 
{\bf Proof} For each  $A \in  \text{SL}(Y, \alpha))$,  the function $f$ itself gives a balanced translation covering  $f: A\circ(X, \omega) \to (\, A\circ(Y, \alpha)\,;p)$.  There is some $\phi \in \text{Aff}^{+}(\,(Y, \alpha)\,)$ giving an isomorphism of translation surfaces of $A\circ(Y,\alpha)$ with $(Y, \alpha)$.  We now proceed as in the previous proof.
\qed

\bigskip

\begin{prop}  \label{ptGivesCover}  With notation as above,  let $m =  [\text{SL}(Y, \alpha; p) :\text{SL}(X, \omega)]$. There is a one-to-$m$ correspondence between elements of the $\text{Aff}^{+}(\,(Y, \alpha)\,)$-orbit of $p$ and balanced coverings $h:  (Z, \beta) \to (Y, \alpha;q)$ such that  $\pi(\, (Z, \beta)\,) \in \mathcal G(X, \omega)$.    Moreover, if 
$(Z, \beta)$ and $(W, \gamma)$ correspond to $q \in \text{Aff}^{+}(\,(Y, \alpha)\,)\cdot p$, then $\pi(\, (Z, \beta)\,) = \pi(\, (W, \gamma) \,)$.  
\end{prop}  

\noindent 
{\bf Proof}  Given $q$ in the $\text{Aff}^{+}(\,(Y, \alpha)\,)$-orbit of $p$, let $\phi(p)=q$.
By Lemma \ref{everyAffHasBldCover}, we find a corresponding $(Z, \beta)$ and $h = \phi\circ f$.   

Now, $\phi$ is unique modulo the stabilizer of  $p$.   Replacing $\phi$ by $\phi\circ \psi$,  with $\psi \in \text{Aff}^{+}(\,(Y, \alpha;p)\,)$ of derivative $B^{-1}$,  upon proceeding as in the proof of Lemma  \ref{everyAffHasBldCover},  we find a balanced covering of translation surfaces $(\phi\circ \psi)\circ f: (W, \gamma)\to (Y, \alpha; q)$, 
with $ (W, \gamma)$ such that  $\pi(\,  (W, \gamma)\,) = \pi(\,AB \circ  (X, \omega)\,)\,.$ Of course,  if $B \in \text{SL}(X, \omega)$, then $f: B\circ(X, \omega) \to B \circ(Y,\alpha;p)$ reduces to the original mapping $f$; thus after the further composition with $\phi$, we then recover the same balanced covering as with $\phi$ by itself.     However, if $B$ is a nontrivial coset representative for $\text{SL}(X, \omega)$ in $\text{SL}(Y, \alpha; p)$, then the reduction to the original mapping cannot be achieved:   $B\circ(X, \omega)$ gives a distinct balanced covering of  $ (Y, \alpha; p)$.   Since $\phi$ is invertible, 
 $(\phi\circ \psi)\circ f$ pull-backs a translation structure from $(Y,\alpha;q)$ that is distinct from that pulled-back by  $\phi \circ  f$.        In any case, each balanced covering $h:  (Z, \beta) \to (Y, \alpha;q)\,$  announces a $(Z, \beta)$  in the $\text{SL}(2, \mathbb R)$-orbit of the translation structure  pulled-back by $\phi\circ f$. Furthermore, we certainly have that $\pi(\, (Z, \beta)\,) \in \mathcal G(X, \omega)$.  

 Now recall that the number of balanced coverings of $(Y, \alpha;q)\,$ is finite and $\text{Aff}^{+}(\,(Y, \alpha;q)\,)$ acts so as to permute them, see say \cite{Vo}.  But,  the pulled-back complex structures in such an orbit are the same:   $\pi(\, (Z, \beta)\,) = \pi(\, (W, \gamma) \,)$.    
\qed

\subsection{Factoring through a Universal Curve}\label{markMod}

\bigskip 
 The (coarse) moduli space of genus $g$ curves with $r$ marked points is denoted by $\mathcal M_{g,r}$.   The forgetful map, given by erasing the final marked point, is  $\pi_{g,r+1}: \mathcal M_{g,r+1} \to \mathcal M_{g,r}$; this is also known as the {\em universal curve}.  See, say, \cite{Hu} and \cite{HM} for these notions.    If $(Z,\beta)$ is of genus $g$ and has $r$ singularities, then since the $\text{SL}(2, \mathbb R)$-action preserves the singularities of $(Z,\beta)$,   the lift of  $\mathcal G(Z,\beta)$ to $\mathcal M_{g,r}$ is a biholomorphic equivalence.    Given a balanced covering $f: (X, \omega) \to (Y, \alpha; \text{Br}(f)\,)$  there is a relationship between $\mathcal G(X, \omega)$ and $\mathcal G(Y, \alpha; \text{Br}(f)\,)$; when $ \text{Br}(f)$ is the union of some singularities and a singleton, then the relationship can be expressed in terms of the universal curve over $\mathcal G(Y, \alpha)$.    
 
  One could discuss the current setting in terms of the generalized Hurwitz spaces as presented in Wewer's dissertation \cite{W}.  
 We give an elementary discussion here; one may turn to \cite{L} for a deeper discussion of correspondences between moduli spaces induced by ramified coverings.   
  We suppose that the genus $h$ translation surface $(Y, \alpha)$ has $r$ singularities, and consider the universal curve over the lift of $\mathcal G(Y, \alpha)$ to $\mathcal M_{h,r}$; the roving $r+1$-st point accounts for the non-singular branch point in coverings related to the initial ramified covering of  $Y$ by $X$.

\begin{thm}\label{fiberUnivCurve}   
 Let $f: (X, \omega)\to (Y, \alpha; p)$ be a balanced covering of translation surfaces as above. 
 Then $\pi: \text{SL}(2, \mathbb R)\circ(X,\omega)\to\mathcal G(X, \omega)$ factors through  the universal curve over $\mathcal G(Y, \alpha)$. 
 \end{thm}  

\noindent 
{\bf Proof}  For general  $A \in \text{SL}(2, \mathbb R)$,    we can replace $f:(X, \omega)\to (Y, \alpha;p)$ by 
$f: A\circ (X, \omega)\to A\circ(Y, \alpha;p)$, and apply Proposition \ref{ptGivesCover} to find the $q$ in the  $\text{Aff}^{+}(\, A\circ(Y,\alpha)\,)$-orbit of $p$ that is  uniquely associated to  $A\circ (X, \omega)$.     Now, since 
$\text{SL}(A\circ(Y,\alpha)\,)$ is conjugated by $A$ to $\text{SL}(Y,\alpha)$, it follows that 
when  $A, B$ are equivalent modulo $\text{SL}(Y, \alpha)$ then the same point of the underlying surface is associated to both $A\circ (X, \omega)$ and $B\circ (X, \omega)$. 
Of course, for each   $ A\circ(Y,\alpha)$, the set of these points $q$ is dense on the complement to the $r$ singularies.  
 Hence, we have a map from the $\text{SL}(2, \mathbb R)$-orbit of   $(X, \omega)$ to a dense subset, $\mathcal C^{0}$,  of $\mathcal C$ the universal curve restricted to 
$\mathcal G(Y, \alpha)$, see Figure 3.

Now, given any of these points $q$,  we can form a balanced covering of the form $h: (Z, \beta) \to (\, A\circ(Y, \alpha); q\,)$, with $\pi(\, (Z, \beta)\,) \in \mathcal G(Y, \alpha)$. 
The composition of these maps does give  $\pi: \text{SL}(2, \mathbb R)\circ(X,\omega)\to\mathcal G(X, \omega)$.\qed  
 \bigskip 

\xymatrix{
& {\text{SL}}(2, \mathbb R)\circ(X,\omega) \ar[dl]\ar[dr] \ar[dd]_{\pi}   \\
 \text{PSL}(X, \omega) \backslash \mathbb H^2 \ar[dr]& & { {\mathcal C}^{0}\,}  \ar[dl] \ar@{^ {(}->}[r] &  {\mathcal C}\, \ar[d] \ar@{^ {(}->}[r] & \mathcal M_{h,r+1}\ar[d]\\
&\mathcal G(X, \omega)&   &\mathcal G(Y,\alpha)\,\ar@{=}[d] \ar@{^ {(}->}[r]        &\mathcal M_{h,r}\ar[d]\\ 
&  &   &\mathcal G(Y,\alpha)\,\ar@{^ {(}->}[r]         & \mathcal M_{h}\\ 
}

\medskip 
\centerline{Figure 3.   Factoring through a cover of the universal curve over $\mathcal G(Y, \alpha)$}

 \bigskip

\begin{prop}  \label{endsAlign}  With notation as above,  suppose that $\mathcal E$ is an end  of $\text{PSL}(X, \omega) \backslash \mathbb H^2$.  Then the image of 
$\mathcal E$ in $\mathcal G(X, \omega)$ lies within the image of a subset 
$\tilde {\mathcal E}$ of $\mathcal C^{0}$ fibering above a cusp of 
$\mathcal G(Y,\alpha)$.    More precisely, if $\mathcal E$ is a cusp, then 
$\tilde {\mathcal E}$ is homeomorphic to $\mathbb R^{+} \times S^1$.   If $\mathcal E$ is an infinite end, then $\tilde {\mathcal E}$ is homeomorphic to $\mathbb R^{+} \times S^1\times S$, with $S$ a dense subset of $S^1$.
\end{prop}  

\noindent 
{\bf Proof}  Suppose first that $\mathcal E$ is a cusp.  Let $\theta = 1/\xi$ be the direction on $(X,\omega)$ corresponding to $\mathcal E$ as in Lemma \ref{fatPtsAreInfEndDir} of \S \ref{endHypSurf}; normalize so that this is the vertical direction.   Of course, $\theta$ is a parabolic direction:  $(X, \omega)$ admits a cylinder composition in this direction, with all cylinders having commensurable inverse moduli.    But, the balanced covering $f: (X, \omega) \to (Y, \alpha; p)$  sends cylinders to cylinders.  Hence,  $\theta$ is also a parabolic direction for  $(Y, \alpha)$, and is such that $p$ splits no cylinder irrationally.   That is, either $p$ lies on the boundary of some cylinder of the decomposition of  $(Y, \alpha)$ in this direction, or else splits some cylinder in a rational manner.    If $p$ lies on the boundary of cylinders, then there is equality of sets of moduli of cylinders of the $(X, \omega)$ and of moduli of cylinders of the $(Y, \alpha)$-decomposition in this direction.   Hence, the minimum value, say $s'$, of $s \in \mathbb R^{+}$ such that $\begin{pmatrix} 1&0\\ s&1\end{pmatrix}$ belongs to  $\text{SL}(Y, \alpha)$ is also the minimum value such that 
this is an element of $\text{SL}(X, \omega)$.   Again since $p$ lies on the boundary of $(Y, \alpha)$-cylinders, the action of each $h_s = \begin{pmatrix} 1&0\\ s&1\end{pmatrix}$ fixes $p$ as a point on the underlying topological surface.   Hence, in accordance with the proof of Theorem \ref{fiberUnivCurve}, $f: h_s \circ (X, \omega) \to h_s \circ (Y, \alpha; p)$ with $s \in [0, s']$ describes a circle in $\mathcal C^{0}$.    There is some $T'>0$ such that the horodisk of index $T'$ projects to the cuspidal region of $\mathcal E$; similarly, there is some $T''$ such that the horodisk of index $T''$ for $(Y, \alpha)$ projects to the cuspidal region of the cusp of $\mathcal G(Y, \alpha)$ corresponding to  the parabolic element $h_s'$.  Let $T = \max\{T', T''\}$.   We find that the image in $\mathcal G(X, \omega)$ of the horodisk of index $T$ of $\mathcal E$ meets the image of the region  $\tilde {\mathcal E} = \{ (\,p, \,(h_s g_t)\circ (Y, \alpha)\,)\,|\, s \in [0, s']\,,\, t >T\,\} \subset \mathcal C^{0}$.  We indeed have that $\tilde {\mathcal E}$ is homeomorphic to $\mathbb R^{+} \times S^1$.

  If $p$ rationally splits some cylinder of $(Y, \alpha)$, then it is still the case that the  single commensurability class of the inverse moduli of the cylinder decompositions of $(X, \omega)$ and of $(Y, \alpha)$ agree.   However,  the parabolic element $h_{s'}$ of 
$\text{SL}(Y, \alpha)$ as above is not an element of   $\text{SL}(X, \omega)$; there is some finite power of this element, thus some $h_{n s'}$ with $n\in \mathbb N$ minimal, that does belong to  $\text{SL}(X, \omega)$.     We find that the horocycle $h_s \circ(X, \omega)$ gives an $n$-fold cover of the corresponding horcycle of $(Y,\alpha)$.   We again can let 
$\tilde {\mathcal E} = \{ (\,p, \,(h_s g_t)\circ (Y, \alpha)\,)\,|\, s \in [0, s']\,,\, t >T\,\} \subset \mathcal C^{0}$.  

Finally, suppose that  $\mathcal E$ is an infinite end.   We continue with the normalization that the vertical direction on $(X, \omega)$ corresponds to $\mathcal E$.    We have a decomposition into cylinders for $(X, \omega)$, which projects to a decomposition of cylinders for $(Y, \alpha; p)$.   There are at least two commensurability classes of inverse moduli of cylinders for the $(X, \omega)$-decomposition (for otherwise $\mathcal E$ would be a cusp) ---  for no $s>0$ is $h_s \in \text{SL}(X, \omega)$.   
Now,  $(Y,\alpha)$ is a Veech surface, hence this  direction of a cylinder decomposition must be parabolic; there is thus a sole commensurability class of inverse moduli for this cylinder decomposition, and we can define $s'$ as above.  Since the balanced covering of $(Y, \alpha; p)$ has at least two such classes, $p$ must irrationally split some cylinder.   But, this introduces at most one new commensurability  class; we find that  $(X, \omega)$ has exactly two commensurability classes of inverse moduli of cylinders for its cylinder decomposition.   
Of more importance,  the orbit of $p$ under powers of $h_{s'}$ is dense in the vertical circle of $(Y, \alpha)$ upon which $p$ lies.    With $T$ as above, we now let $\tilde {\mathcal E} = \{ (\,h_s(p), \,(h_s g_t)\circ (Y, \alpha)\,)\,|\, s \in [0, s']\,,\, t >T\,\} \subset \mathcal C^{0}$.    This is is homeomorphic to $\mathbb R^{+} \times S^1\times S$, with $S$ a dense subset of $S^1$.
\qed

\bigskip 
Theorem \ref{thmBigCusp} now follows --- each end of $\mathcal G(X, \omega)$ must be covered by some end of $\text{PSL}(X, \omega)\backslash \mathbb H^2$, but each of these is given by one of finitely many ends:  the image of the fibering over one of the finitely many ends of $\mathcal G(Y, \alpha)$.

\end{document}